\documentclass[a4paper,11pt]{amsart}
\usepackage{amscd,amssymb,bbm}
\usepackage[all]{xy}

\newcommand{\B}{\mathcal B}

\renewcommand{\L}{\Lambda}
\newcommand{\Z}{\mathbb Z}

\newcommand{\ainf}{A_{\infty}}
\newcommand{\Aut}{\operatorname{Aut}}

\newcommand{\cone}{\operatorname{cone}}
\newcommand{\comj}{\mathsf{comproj}}
\newcommand{\Db}{\mathcal D ^b}
\newcommand{\dege}{\leq_{\operatorname{deg}}}
\newcommand{\empt}{\varnothing}
\newcommand{\Ext}{\operatorname{Ext}}
\newcommand{\Gl}{\operatorname{\mathsf{GL}}}
\newcommand{\Hom}{\operatorname{Hom}}
\renewcommand{\Im}{\operatorname{Im}}
\renewcommand{\mod}{\operatorname{mod}}
\newcommand{\Mod}{\operatorname{Mod}}
\newcommand{\op}{{\operatorname{op}}}

\newcommand{\pil}{\rightarrow}

\newcommand{\Rep}{\operatorname{Rep}}
\renewcommand{\top}{\operatorname{top}}
\newcommand{\xpil}{\xrightarrow}

\newtheorem{lem}{Lemma}[section]
\newtheorem{prop}[lem]{Proposition}
\newtheorem{cor}[lem]{Corollary}
\newtheorem{thm}[lem]{Theorem}
\newtheorem*{thm2}{Theorem}
\theoremstyle{definition}

\newtheorem{remark}[lem]{Remark}

\title{Degeneration of A-infinity modules}
\author{Bernt Tore Jensen}
\address{Institutt for matematiske fag\\ NTNU\\ NO--7491 Trondheim\\ Norway}
\email{berntj@math.ntnu.no}

\author{Dag Madsen}
\address{Institutt for matematiske fag\\ NTNU\\ NO--7491 Trondheim\\ Norway}
\email{dagma@math.ntnu.no}

\author{Xiuping Su}
\address{Mathematisches Institut\\ Universit\"{a}t zu K\"{o}ln\\
Weyertal 86-90\\ 50931 K\"{o}ln\\ Germany}
\email{xsu@math.uni-koeln.de}

\subjclass[2000]{Primary 18E30; Secondary 14L30, 16G10}

\keywords{A-infinity modules, Degeneration (Orbit closure), Derived
categories}

\begin{document}

\begin{abstract}
In this paper we use $\ainf$-modules to study the derived category
of a finite dimensional algebra over an algebraically closed field.
We study varieties parameterising $\ainf$-modules. These varieties
carry an action of an algebraic group such that orbits correspond to
quasi-isomorphism classes of complexes in the derived category. We
describe orbit closures in these varieties, generalising a result of
Zwara and Riedtmann for modules.
\end{abstract}

\maketitle

\section*{Introduction}

In this paper we study parameter spaces of objects in the derived
category of a finite dimensional algebra $\Lambda$ over an algebraically
closed field $k$.

The derived category of $\L$ can be defined by formally inverting
all quasi-isomorphisms in the category of complexes of $\L$-modules.
Recall that a quasi-isomorphism is a chain map which induces
isomorphism in homology. So two isomorphic objects in the derived
category have isomorphic homologies, but the converse is in general
not true. A natural point of view is therefore that
quasi-isomorphism classes are determined by the isomorphism type of
their homology and some extra information which explains how the
homology is put together to form a complex. In this paper we will
take this point of view.

We study families of objects in the derived category with fixed
dimensions or even fixed isomorphism type in homology. We do so by
replacing the ordinary derived category with the equivalent derived
category of $\ainf$-modules \cite{Ainf}. The equivalence functor
preserves homology, so dimensions and isomorphism types are not
affected. The $\ainf$-modules are called polydules in \cite{Doc},
and we use both names in this paper.

Keller has observed \cite{Ainf} that by fixing dimensions in
homology $\underline{d}=(d_a, \ldots,d_b)$ we obtain a variety, here
denoted by $\mathcal{M}=\mathcal{M}_{\Lambda,{\underline d}}$, which
parameterises polydules with homology dimensions $\underline d$.
There is an algebraic group $\mathcal{G}=\mathcal{G}_{\underline d}$
acting on $\mathcal{M}_{\L,{\underline d}}$ such that orbits
correspond to quasi-isomorphism classes of polydules. It is this
variety and group action we examine in the present paper.

The main result is a generalisation of a result of Zwara \cite{Zwa}
and Riedtmann \cite{Riedt} from module varieties to these varieties
of polydules.

\begin{thm2} Let $M,N\in \mathcal M_{{\L, \underline d}}$. Then $N\in
\overline{\mathcal{G}\cdot M}$ if and only if there is a finite
dimensional $\L$-polydule $Z$ and an exact triangle $$Z \pil M
\oplus Z \pil N \pil Z[1].$$
\end{thm2}

The paper is organised as follows. In section 1 we recall the
definition of $\ainf$-modules and their morphisms. In section 2 we
give a summary of the most important known results concerning
$\ainf$-modules and their derived category. In section 3 we describe
in detail the variety $\mathcal{M}$ and the group $\mathcal{G}$ with
its action on $\mathcal M$. We consider the special cases where the
homology is given as representations of a quiver with relations or
where the homology module structure is fixed in section 4. In
section 5 we prove our main result, an algebraic characterisation of
degeneration. We give several examples in section 6.

\section{A-infinity modules}\label{defi}

In this section we give the basic definitions of $\ainf$-modules and
their morphisms. We stress that our definitions are adapted to the
case we are interested in, namely $\ainf$-modules over associative
algebras, and are only valid for this case. For the definitions in
the most general case ($\ainf$-modules over $\ainf$-algebras), we
refer to \cite{Ainf} and \cite{Doc}.

For our purposes we find it most practical to use a left module
notation. A (left) $\ainf$-module over an algebra $\L$ in our sense
is the same as a (right) $\ainf$-module over $\L^{\op}$ in the sense
of \cite{Ainf}. Since all our algebra elements are in degree $0$ we
avoid certain sign issues.

\subsection{Modules}

Let $k$ be a commutative ring and $\L$ an associative $k$-algebra.
An \emph{$\ainf$-module} over $\L$ (or a \emph{$\L$-polydule}
\cite{Doc}) is a graded $k$-module
$$M=\bigoplus_{i \in \Z} M^i$$ with $k$-linear graded maps
$$m_n \colon \L^{(\otimes_k) n-1} \otimes_k M \pil M,\textit{ } n \geq 1$$
of degree $n-2$ satisfying the rules
\begin{align*}
m_1 m_1 &=0,\\
m_1 m_2 &=m_2(\mathtt 1 \otimes m_1),\\
\intertext{and for $n \geq 3$}
\sum_{i=1}^n(-1)^{i(n-1)}m_{n-i+1}(\mathtt 1^{\otimes n-i} \otimes
m_i) &=\sum_{j=1}^{n-2}(-1)^{j-1} m_{n-1}(\mathtt 1^{\otimes n-j-2}
\otimes \mu \otimes \mathtt 1^{\otimes j}),
\end{align*}
where $\mathtt 1$ is the identity map and $\mu$ is the
multiplication in $\L$.

In particular $(M,m_1)$ is a complex of $k$-modules. Note that if
$M$ is bounded, then $m_n=0$ for $n \gg 0$.

\subsection{Change of base rings}

Suppose $M=\bigoplus_{i \in \Z} M^i$ a $\L$-polydule with maps
$m_n^M$, $n \geq 1$. If $V$ is an (ungraded) $k$-module, then the
graded $k$-module $$M \otimes_k V=\bigoplus_{i \in \Z} (M^i
\otimes_k V)$$ is a $\L$-polydule with maps $$m_n^{M \otimes V}=
m_n^M \otimes \mathtt 1_V \colon \L^{(\otimes_k) n-1} \otimes_k M
\otimes_k V \pil M \otimes_k V,\textit{ } n \geq 1.$$

Suppose now that $\theta \colon k \pil k'$ is a morphism of
commutative rings. Then $\L \otimes_k k'$ is a $k'$-algebra, and a
$\L \otimes_k k'$-polydule is a graded $k'$-module $N=\bigoplus_{i
\in \Z} N^i$ with $k'$-linear maps $$m_n^N \colon (\L \otimes_k
k')^{(\otimes_{k'}) n-1} \otimes_{k'} N \pil N ,\textit{ } n \geq
1$$ satisfying the equations. Since $$(\L \otimes_k
k')^{(\otimes_{k'}) n-1} \otimes_{k'} N \simeq \L^{(\otimes_{k})
n-1} \otimes_k N,$$ this structure is determined by $k'$-linear maps
$$m_n^N \colon \L^{(\otimes_{k}) n-1} \otimes_k N \pil N,\textit{ } n \geq
1.$$ The map $\theta \colon k \pil k'$ induces a $\L$-polydule
structure on $N$.

If $M=\bigoplus_{i \in \Z} M^i$ a $\L$-polydule, there is a $\L
\otimes_k k'$-polydule structure on the graded $k'$-module
$$M \otimes_k k'=\bigoplus_{i \in \Z} (M^i \otimes_k k')$$ given by
the $k'$-linear maps $$m_n^{M \otimes k'}=m_n^M \otimes \mathtt 1
\colon \L^{(\otimes_k) n-1} \otimes_k M \otimes_k k' \pil M
\otimes_k k',\textit{ } n \geq 1.$$

\subsection{Morphisms and compositions}

A \emph{morphism} $f \colon L \pil N$ between two $\L$-polydules $L$
and $N$ is given by a family of maps
$$f_n \colon \L^{(\otimes_k) n-1} \otimes_k L \pil N, \textit{ } n \geq
1$$ of degree $n-1$ satisfying the rules
\begin{align*}
f_1 m_1 &=m_1 f_1,\\
f_1 m_2-f_2(\mathtt 1 \otimes m_1) &=m_2(\mathtt 1 \otimes f_1)+m_1
f_2,\\
\intertext{and for $n \geq 3$}
\sum_{i=1}^n(-1)^{i(n-1)}f_{n-i+1}(\mathtt 1^{\otimes n-i} \otimes
m_i) &+\sum_{j=1}^{n-2}(-1)^j f_{n-1}(\mathtt 1^{\otimes n-j-2}
\otimes \mu \otimes \mathtt 1^{\otimes j})\\ &=\sum_{r=1}^n
(-1)^{(r+1)n} m_{n-r+1}(\mathtt 1^{\otimes n-r} \otimes f_r).
\end{align*}
Note in particular that $f_1$ is a chain map $f_1 \colon (L,m_1)
\pil (N,m_1)$ between complexes of $k$-modules. A morphism is called
a \emph{quasi-isomorphism} if $f_1 \colon (L,m_1) \pil (N,m_1)$ is a
quasi-isomorphism of complexes.

A morphism $f$ is called \emph{strict} if $f_i=0$ for $i>1$.
Furthermore $f$ is called a \emph{strict monomorphism} if $f$ is
strict and $f_1$ is a monomorphism. Similarly $f$ is called a
\emph{strict epimorphism} if $f$ is strict and $f_1$ is an
epimorphism.

The \emph{identity morphism} $\mathtt 1_N \colon N \pil N$ is given
by $f_1=\mathtt 1$ and $f_i=0$ for all $i>0$. The \emph{composition}
$fg \colon N \pil M$ of two morphisms $f \colon L \pil M$ and $g
\colon N \pil L$ is given by the rule
$$(fg)_n=\sum_{i=1}^n (-1)^{(i-1)n}f_{n-i+1}(\mathtt 1^{\otimes n-i} \otimes
g_i).$$

A morphism $f \colon L \pil N$ is called \emph{null-homotopic} if
there is a family of maps $$h_n \colon \L^{(\otimes_k) n-1}
\otimes_k L \pil N, \textit{ } n \geq 1$$ of degree $-n$ such that
\begin{align*}
f_1= &m_1 h_1+h_1 m_1,\\
f_2= &m_1 h_2-m_2(\mathtt 1 \otimes h_1)-h_2(\mathtt 1 \otimes
m_1)+h_1 m_2\\
\intertext{and for $n \geq 3$} f_n= &\sum_{r=1}^n(-1)^{r(n-1)}
m_{n-r+1}(\mathtt 1^{\otimes n-r} \otimes
h_r)\\&+\sum_{i=1}^n(-1)^{i(n-1)}h_{n-i+1}(\mathtt 1^{\otimes n-i}
\otimes m_i)\\&+\sum_{j=1}^{n-2}(-1)^j h_{n-1}(\mathtt 1^{\otimes
n-j-2} \otimes \mu \otimes \mathtt 1^{\otimes j}).
\end{align*}

To morphisms $f$ and $g$ from $L$ to $N$ are called \emph{homotopic}
if $f-g$ is null-homotopic. In this case we say that $h$ is an
\emph{homotopy} between $f$ and $g$. A morphism $f \colon L \pil N$
is called a \emph{homotopy equivalence} if there exists a morphism
$g \colon N \pil L$ such that $gf$ is homotopic to $\mathtt 1_L$ and
$fg$ is homotopic to $\mathtt 1_N$.

\subsection{Unitality}

A $\L$-polydule $M$ is called \emph{strictly unital} if for all $a
\in M$ $m_2(1,a)=a$ and $m_n(\lambda_1, \ldots,\lambda_{n-1},a)=0$
if $n \geq 3$ and $1 \in \{\lambda_1,\ldots,\lambda_{n-1} \}$.

A morphism $f \colon L \pil N$ is called \emph{strictly unital} if
$f_n(\lambda_1, \ldots,\lambda_{n-1},a)=0$ whenever $n \geq 2$ and
$1 \in \{\lambda_1,\ldots,\lambda_{n-1} \}$. Similarly an homotopy
$h$ between two morphisms $f$ and $g$ is called \emph{strictly
unital} if $h_n(\lambda_1, \ldots,\lambda_{n-1},a)=0$ whenever $n
\geq 2$ and $1 \in \{\lambda_1,\ldots,\lambda_{n-1} \}$. If there is
such a homotopy between $f$ and $g$, we denote this fact by $f
\sim_h g$.

\section{The derived category}

Suppose now that $k$ is a field and $\L$ a $k$-algebra. In this
section we give a summary of the, for our purposes, most important
results from \cite{Ainf} and \cite{Doc} concerning $\ainf$-modules
and their derived category.

Let $\Mod_\infty \L$ denote the category of strictly unital
$\L$-polydules and strictly unital morphisms. The \emph{homotopy
category} $\mathcal H_\infty(\L)$ is defined to be the quotient
category $ \mathcal H_\infty(\L)=\Mod_\infty \L / \sim_h$, so
morphisms in $\mathcal H_\infty(\L)$ are homotopy classes of
morphisms in $\Mod_\infty \L$. The following theorem is of central
importance.

\begin{thm} \cite[4.2]{Ainf} \cite[2.4.1.1]{Doc}
A morphism $f \colon L \pil N$ is a homotopy equivalence if and only
if it is a quasi-isomorphism.
\end{thm}

This gives us one out of several essentially equivalent ways of
defining the \emph{derived category} $\mathcal D_\infty(\L)$ (see
section 4.1.3 in \cite{Doc}). We simply define it to be the same as
the homotopy category, so $\mathcal D_\infty(\L)=\mathcal
H_\infty(\L)$.

The category $\mathcal D_\infty(\L)$ is a triangulated category. The
\emph{shift} of a $\L$-polydule $M$, denoted $M[1]$, is the
$\L$-polydule defined by $(M[1])^i=M^{i+1}$ and $m_n^{M[1]}=(-1)^n
m_n^M$. If $f \colon L \pil N$ is a morphism in $\mathcal
D_\infty(\L)$, then there exists an object $\cone f$, called the
\emph{cone} of $f$, and an exact triangle $$L \xpil f N \xpil j
\cone f \xpil p L[1]$$ in $\mathcal D_\infty(\L)$. The cone has
graded parts
$$(\cone f)^i=N^i \oplus (L[1])^i,$$ and the $\L$-polydule
structure is defined by $$m_n^{\cone f}=\begin{pmatrix} m_n^N &
(-1)^n f_n\\ 0 & m_n^{M[1]}
\end{pmatrix}.$$ Every exact triangle in $\mathcal D_\infty(\L)$
is isomorphic to one
of the form above. The morphism $j$ is a strict monomorphism, the
morphism $p$ is a strict epimorphism and the sequence of graded
$k$-vector spaces $0 \pil N \xpil {j_1} \cone f \xpil {p_1} L[1]
\pil 0$ is exact.

Using the rotation axiom, we get the same description of exact
triangles as in \cite[5.2]{Ainf}. Let $L \xpil j M \xpil p N$ be a
sequence with $j$ a strict monomorphism, with $p$ a strict
epimorphism and where the sequence of graded $k$-vector spaces $0
\pil L \xpil {j_1} M \xpil {p_1} N \pil 0$ is exact. We call such a
sequence an \emph{exact sequence of $\L$-polydules}. In this case
there is an exact triangle $L \xpil j M \xpil p N \pil L[1]$ in
$\mathcal D_\infty(\L)$, so $\cone j \simeq N$ in $\mathcal
D_\infty(\L)$. All exact triangles are obtained in this way.

If $(C,\delta)$ is a complex of $\L$-modules, then $C$ can be viewed
as a $\L$-polydule where $m_1=\delta$, the map $m_2 \colon \L
\otimes C \pil C$ is given by the $\L$-module structure and $m_n=0$
for $n \geq 3$. Identifying complexes with polydules in this way
induces a triangle functor $\mathcal D(\Mod \L) \pil \mathcal
D_\infty(\L)$, where $\mathcal D(\Mod \L)$ is the usual (unbounded)
derived category of $\L$. It turns out that this functor is an
equivalence.

\begin{thm} \cite[4.3]{Ainf} \cite[2.4.2.3]{Doc}
The canonical functor $\mathcal D(\Mod \L) \pil \mathcal
D_\infty(\L)$ is an equivalence of triangulated categories.
\end{thm}

In other words each quasi-isomorphism class of $\L$-polydules
corresponds to (has as a subclass) exactly one quasi-isomorphism
class of complexes of $\L$-modules. So every polydule is
quasi-isomorphic to a complex in this sense.

In this paper however, we will focus on another type of polydules,
the ones with $m_1=0$. A good thing about such objects is that
quasi-isomorphisms between them are invertible already in
$\Mod_\infty \L$.

\begin{thm} \cite[4.3]{Ainf}
Suppose $M$ and $N$ are two $\L$-polydules with $m_1^M=0=m_1^N$.
Then $M$ and $N$ are quasi-isomorphic if and only they are
isomorphic.
\end{thm}

Each quasi-isomorphism class of polydules has an object with
$m_1=0$, so polydules of this type can be chosen as representatives
of isomorphism classes in the derived category $\mathcal
D_\infty(\L)$.

\begin{thm}\cite[3.3.1.7]{Doc}
Suppose $M$ is a strictly unital $\L$-polydule. Then there is a
strictly unital $\ainf$-module $M'$ with $m_1^{M'}=0$ and a strictly
unital quasi-isomorphism $M' \pil M$.
\end{thm}

If $M$ is a complex, meaning $m_n^M=0$ for $n \geq 3$, then the
graded parts of $M'$ are the homology modules of $(M,m_1^M)$. In
general if $N$ is a $\L$-polydule with $m_1^N=0$, then each $N^i$,
$i \in \Z$, is a $\L$-module in its own right where the module
structure is $m_2^N$ restricted to $m_2^N \colon \L \oplus N^i \pil
N^i$. It is also worth noting that if $f \colon M \pil N$ is a
morphism and $m_1^M=0=m_1^N$, then for each $i \in \Z$, the map $f_1
\colon M^i \pil N^i$ is a $\L$-module morphism. A polydule with
$m_1=0$ can informally be thought of as "homology + structure".

Our aim is to study the bounded derived category of finitely
generated $\L$-modules, $\Db (\mod \L)$, when $\L$ is a finite
dimensional $k$-algebra. The category $\Db (\mod\L)$ is in this case
equivalent to the full subcategory of $\mathcal D (\Mod \L)$
consisting of objects with finite dimensional total homology. We get
the following characterisation of the corresponding polydules.

\begin{prop}
Let $\L$ be a finite dimensional $k$-algebra. A $\L$-polydule $M$ is
in the essential image of the composed functor $\Db (\mod \L)
\hookrightarrow \mathcal D(\Mod \L) \xpil \sim \mathcal
D_\infty(\L)$ if and only if it is quasi-isomorphic to a
$\L$-polydule $M'$ with $m_1^{M'}=0$ and $\dim_k M' < \infty$.
\end{prop}

We denote by $\mathcal D^f_\infty(\L)$ the full subcategory of
$\mathcal D_\infty(\L)$ consisting of objects in this image.
Obviously we have an equivalence of triangulated categories $\Db
(\mod \L) \xpil \sim \mathcal D^f_\infty(\L)$.

\section{Description of variety}

In this section we describe the geometric object we want to study.
Suppose that $k$ is an algebraically closed field and $\L$ is a
finite dimensional $k$-algebra. We keep these assumptions for the
rest of the paper. We fix dimensions in homology, non-zero only on a
finite interval, and make a variety out of the possible
$\L$-polydule structures with $m_1=0$ and the given homology
dimensions. Via the equivalence $\Db (\mod \L) \xpil \sim \mathcal
D^f_\infty(\L)$, we can view this as a parameter space of objects in
the derived category with the given homology dimensions. There is a
group of quasi-isomorphisms acting on this space. As noted by Keller
\cite[4.3]{Ainf}, what we have in this situation is an algebraic
group acting on an algebraic variety. While this section is devoted
to definitions, in section \ref{orbit} we pursue our main goal,
namely to find an algebraic characterisation of orbit closure
(degeneration).

Suppose $\B=\{ \nu_1=1,\nu_2, \ldots, \nu_t \}$ is a $k$-basis of
$\L$. We fix integers $a < b$ and a finite dimensional graded
$k$-vector space $V=\oplus_{i=a}^b V^i$. We want to construct a
variety $\mathcal M$ of possible $\ainf$-module structures on $V$.
Let $d_i=\dim_k V^i$ for $a \leq i \leq b$, and let $\underline d$
denote the vector $\underline d=(d_a, \ldots, d_b)$. We fix bases
for the graded parts of $V$, so that maps between them are
represented by matrices.

A point in the variety $\mathcal M=\mathcal M_{\L,\underline d}$ is
a collection of matrices $M(S,i)$, one for each sequence
$S=(\lambda_{|S|}, \ldots, \lambda_2, \lambda_1)$ of $1 \leq |S|
\leq b-a+1$ elements from $\B$ and $i$ a number satisfying $a+|S|-1
\leq i \leq b$. The matrix $M{(S,i)}$ is of size $d_{i-|S|+1} \times
d_i$ and represents the map
$$m_{|S|+1}(\lambda_{|S|}, \ldots, \lambda_1,-) \colon V^i \pil V^{i-|S|+1}.$$
We consider only strictly unital polydules, so if $1 \in S$ then
$$M(S,i)=\left \{
\begin{array}{c c}
I & \text { if } |S|=1\\
0 & \text { if } |S|>1
\end{array} \right . $$

Given a sequence $S=(\lambda_{|S|}, \ldots,\lambda_1)$ and $1 \leq l
\leq |S|-1$, there are uniquely determined elements $c^{S,l}_r \in
k$, $1 \leq r \leq t$ satisfying the expression
$$\lambda_{l+1} \lambda_l=\sum_{r=1}^t c^{S,l}_r \nu_r.$$
The map $$m_{|S|}(\lambda_{|S|},\ldots,\lambda_{l+2},\lambda_{l+1}
\lambda_l,\lambda_{l-1},\ldots,\lambda_1,-) \colon V^i \pil
V^{i-|S|+2}$$ is represented by the matrix $\sum_{r=1}^t c^{S,l}_r
M(S_{l/r},i)$, where $S_{l/r}$ is the sequence
$$S_{l/r}=(\lambda_{|S|},
\ldots,\lambda_{l+2},\nu_r,\lambda_{l-1},\ldots,\lambda_1)$$ of
length $|S|-1$ obtained from $S$ by deleting $\lambda_{l+1}$ and
$\lambda_l$ and inserting $\nu_r$ in position $l$ from the right.

The matrices $M(S,i)$ have to satisfy the relations
$$\sum_{S=[S',S'']} (-1)^{|S|(|S''|+1)} M{(S',i-|S''|+1)} \cdot
M{(S'',i)} =\sum_{l=1}^{|S|-1} \sum_{r=1}^t (-1)^{l-1} c^{S,l}_r
\cdot M(S_{l/r},i)$$ (This and other formulas in this section are
obtained from the corresponding formulas in section \ref{defi}.),
where the sum on the left hand side is taken over all decompositions
of $S$ into two parts $S'=(\lambda_{|S|},\ldots,\lambda_{|S''|+1})$
and $S''=(\lambda_{|S''|},\ldots,\lambda_1)$ with $|S'|,|S''| \geq
1$. We have one such relation for each pair $(S,i)$ with $2 \leq |S|
\leq b-a+2$ and $a+|S|-2 \leq i \leq b$.

\begin{remark}
Our construction of the variety $\mathcal M_{\L,\underline d}$
depends on the fact that $\L$ has a finite $k$-basis. It might also
be possible to construct a similar variety in some cases when $\L$
is not finite dimensional.
\end{remark}

Next we describe the group $\mathcal G=\mathcal G_{\underline d}$ of
quasi-isomorphisms. An element $g$ in $\mathcal G$ is a collection
of matrices $F(S,i)$, one for each sequence $S$ of $0 \leq |S| \leq
b-a$ elements and $i$ a number satisfying $a+|S| \leq i \leq b$.
Note that here we allow $S$ to be the empty sequence $\empt$. The
matrix $F(S,i)$ is of size $d_{i-|S|} \times d_i$ and is thought of
as representing a map
$$f_{|S|+1}(\lambda_{|S|},\ldots,\lambda_1,-) \colon V^i \pil
V^{i-|S|}.$$ We only want to consider strictly unital morphisms, so
we add the condition that $F(S,i)=0$ whenever $1 \in S$. Since $f_1
\colon V \pil V$ is supposed to be a quasi-isomorphism, the matrices
$F(\empt,i)$ must be invertible for all $a \leq i \leq b$. These are
the only conditions for $g$ to be an element in $\mathcal G$, there
are no other conditions.

The product $g=(F(S,i))$ of two group elements $g'=(F'(S,i))$ and
$g''=(F''(S,i))$ is given by the rule
$$F(S,i)=\sum_{S=[S',S'']}(-1)^{|S'||S''|}F'(S',i-|S|) \cdot
F''(S'',i).$$

The identity element $e$ of this group is the element given by
$F(\empt,i)=I_{d_i}$ for all $i$ and $F(S,i)=0$ for all $S \neq
\empt$. The inverse $g^{-1}=(F'(S,i))$ of an element $g=(F(S,i))$ is
given recursively by
$$F'(\empt,i)=F(\empt,i)^{-1}$$ and
$$F'(S,i)=\Bigl [ \sum_{S=[S',S''],S'' \neq \empt}(-1)^{|S'||S''|+1}
F'(S',i-|S|) \cdot F(S'',i) \Bigl ]
F(\empt,i)^{-1}$$ for all $S \neq \empt$.

We now describe the action of $\mathcal G$ on $\mathcal M$. The map
$\ast \colon \mathcal G \times \mathcal M \pil \mathcal M$ is given
in such a way that if $g=(F(S,i))$ and $x=(M(S,i))$, then $x'=g \ast
x$ is given by
\begin{align*}
M'(S,i)=&\Bigl{[}\sum_{S=[S',S''],S'' \neq \empt} (-1)^{|S'|(|S''|+1)}
F(S',i-|S''|+1) \cdot M(S'',i)\\
&+\sum_{l=1}^{|S|-1} \sum_{r=1}^t (-1)^{l-1} c^{S,l}_r \cdot
F(S_{l/r},i)\\&+\sum_{S=[S',S''],S',S'' \neq \empt}
(-1)^{|S'||S''|+1} M'(S',i-|S|) \cdot F(S'',i)\Bigl]F(\empt,i)^{-1}.
\end{align*}

The $\mathcal G$-orbits in $\mathcal M$ correspond to isomorphism
classes of objects in $\Db(\mod \L)$ with the given dimensions in
homology. We define degeneration in $\mathcal M$ in the usual way as
follows. Let $M$ and $N$ be two $\L$-polydules with $m_1^M=0=m_1^N$
and the given homology dimensions. We use the same letters to denote
the corresponding points in $\mathcal M$. If $N$ belongs to the
Zariski closure of the $\mathcal G$-orbit of $M$, in mathematical
notation $N \in \overline {\mathcal G \cdot M}$, we say that $M$
\emph{degenerates to $N$} and we denote this fact by $M \dege N$.
The relation $\dege$ is a partial order on the set of isomorphism
classes of $\L$-polydules with $m_1=0$ and the given homology
dimensions.

\section{Variations and simplifications}\label{simp}

The variety $\mathcal{M}$ with its action by $\mathcal{G}$ has a
rather complicated definition which is not so easy to work with in
general. In this section we consider some special cases of algebras
and polydules where we can define smaller and simpler spaces and
group actions. We will illustrate these simplifications on some
concrete examples in section \ref{ex}.

Let $\mod(A,\underline{d})\subseteq
\mathcal{M}_{\Lambda,\underline{d}}$ be the subvariety consisting of
polydules with $M(S,i)=0$ for all $S$ with $|S|>1$. That is, each
polydule in $\mod(\Lambda,\underline{d})$ is an graded
$\Lambda$-module, and $\mod(\Lambda,\underline{d})$ is isomorphic to
a product of module varieties
$$\mod(\Lambda,\underline{d})\cong \prod_i \mod(\Lambda,d_i).$$

There is a projection map
$$\pi \colon \mathcal{M}_{\Lambda,\underline{d}}\longrightarrow \mod(\L,\underline{d})$$
defined by $\pi(M(S,i))=0$ for all $S$ with $|S|>1$. Here $M$ and
$N$ belong to the same fibre if and only if they have equal homology
as graded modules. That is, $\pi(M)=H^*M$ computes the homology of a
polydule.

The group $\Gl_{\underline d}=\prod_i \Gl_{d_i}$ acts on
$\mod(\Lambda, \underline{d})$ by conjugation, and the orbits are in
bijection with isomorphism classes of graded $\Lambda$-modules with
graded dimension $\underline{d}$. There is a morphism $$\mathcal{G}
\longrightarrow \Gl_{\underline d}$$ of algebraic groups mapping all
$F(S,i)$ with $|S|>0$ to zero. We let $\mathcal{G}$ act on
$\mod(\Lambda,d)$ via this morphism. Then by the definition of the
action of $\mathcal{G}$ on $\mathcal{M}$, the projection $\pi$ is
$\mathcal{G}$-equivariant.

\subsection{Quivers with relations}\label{quiv}

A very important combinatorial tool for studying
modules over finite dimensional algebras is to use quivers with
relations. We shall see that quivers are also useful when studying
polydules.

Now assume that $\Lambda=kQ/I$ for a quiver $Q=(Q_0,Q_1)$ and an
admissible ideal $I\subseteq kQ$. Then it is well known that the
category of $\Lambda$-modules is equivalent to the category of
representations of $Q$ satisfying the relations in $I$ (see
\cite{Ars} for an introduction to quivers and their
representations). The change of point of view from modules to
representations can be regarded as a change of underlying monoidal
base category from $k$-modules to $(k^{\times r})$-modules, where
$r=|Q_0|$. With the new base category the general theory in
\cite{Doc} is still valid, and a polydule in the new sense is the
same a graded representation of $Q$ with an $A_\infty$-structure
respecting the quiver. We call such a graded representation a
\emph{polyrepresentation}. From the general theory \cite{Doc} it
follows that the derived category of polyrepresentations is triangle
equivalent to the derived category of representations, which we know
is equivalent to the derived category of $\L$-modules.

We will now define a variety of polyrepresentations. We fix not only
dimensions in homology, but dimensions vectors. So the given data is
a sequence $\underline{e}=(\underline{e}^a,\ldots,\underline{e}^b)$
of dimension vectors $\underline e^i=(e^i_1,\ldots,e^i_r)$, $a \leq
i \leq b$. Let $\underline d=(d_a, \ldots, d_b)$, where $d_i=\sum_{j
\in Q_0}e_j^i$. Assume a basis for $\L$ is chosen that consists only
of paths. We define a subvariety
$$\Rep_\infty(\Lambda,\underline{e}) \subseteq \mathcal
M_{\L,\underline d}$$ as follows. For each point $x=(M(S,i)) \in
\mathcal M_{\L,\underline d}$, divide each matrix $M(S,i)$ into
blocks $C_{j,l}(S,i)$, $1 \leq j,l \leq r$, of size $e^{i-|S|+1}_l
\times e^i_j$. The subvariety $\Rep_\infty(\Lambda,\underline{e})$
consists of all $x \in \mathcal M_{\L,\underline d}$ which satisfy
all of the following three conditions.
\begin{itemize}
\item[(i)] $M(S,i)=0$
if the composition $\lambda_{|S|}\cdots\lambda_1$ of the paths in
$S$ is not a path in $Q$.
\item[(ii)] $C_{j,l}(S,i)=0$ if the
composition $\lambda_{|S|}\cdots\lambda_1$ of the paths in $S$ is a
path in $Q$ starting in vertex $s$ and ending in vertex $t$ and
$(j,l) \neq (s,t)$.
\item[(iii)] $M(S,i)=0$ if $|S|>1$ and $S$
contains a path of length zero.
\end{itemize}
The variety $\Rep_\infty(\Lambda,\underline{e})$ is the variety of
polyrepresentations with sequence of dimension vectors
$\underline{e}=(\underline{e}^a,\ldots,\underline{e}^b)$.

The image of $\Rep_\infty(\Lambda,\underline{e})$ under $\pi$,
$$\pi(\Rep_\infty(\L,\underline{e}))=\Rep(\L,\underline{e}) \subseteq
\mod(\Lambda,\underline{d}),$$ is the variety of graded
representations of $Q$ satisfying the relations in $I$ with sequence
of dimension vectors
$\underline{e}=(\underline{e}^a,\ldots,\underline{e}^b)$. For each
$a \leq i \leq b$, a further projection to the $i$th component gives
$$\Rep(\Lambda,\underline{e}^i) \subseteq \mod(\Lambda,d_i),$$ the
variety of representations of $Q$ satisfying the relations in $I$
with dimension vector $\underline e^i$.

For each element $g=(F(S,i)) \in \mathcal G_{\underline d}$, divide
each matrix $F(S,i)$ into blocks $D_{j,l}(S,i)$, $1 \leq j,l \leq
r$, of size $e^{i-|S|}_l \times e^i_j$. Let
$\mathcal{G}_{Q}=\mathcal{G}_{Q,\underline e}$ be the subgroup of
$\mathcal{G}_{\underline d}$ consisting of the elements $g \in
\mathcal{G}_{\underline d}$ which satisfy all of the following four
conditions.
\begin{itemize}
\item[(i)] $D_{j,l}(\empt,i)=0$ if $j \neq l$.
\item[(ii)] $F(S,i)=0$
if $|S|>0$ and the composition $\lambda_{|S|}\cdots\lambda_1$ of the
paths in $S$ is not a path in $Q$.
\item[(iii)] $D_{j,l}(S,i)=0$ if $|S|>0$, the
composition $\lambda_{|S|}\cdots \lambda_1$ of the paths in $S$ is a
path in $Q$ starting in vertex $s$ and ending in vertex $t$ and
$(j,l) \neq (s,t)$.
\item[(iv)] $F(S,i)=0$ if $|S|>0$ and $S$
contains a path of length zero.
\end{itemize}
The group $\mathcal{G}_{Q}$ is the group of quasi-isomorphisms which
respect the quiver. By the general theory \cite{Doc}, we know that
the $\mathcal{G}_{Q}$-orbits in $\Rep_\infty(\Lambda,\underline{e})$
correspond to quasi-isomorphism classes of polyrepresentations and
therefore to isomorphism classes of objects in $\Db(\mod \L)$.
Moreover, $\mathcal{G}\cdot \Rep_\infty(\Lambda,\underline{e})$ is a
connected component of $\mathcal{M}_{\Lambda,\underline d}$, and
each connected component is obtained in this way. This follows,
since each polydule degenerates to its homology, and each module
degenerates to the direct sum of its composition factors. We give a
short argument showing that a polydule degenerates to its homology.

\begin{lem} \label{deglemma}
Let $M \in \mathcal{M}_{\Lambda,\underline{d}}$. Then $M\dege
H^*(M)$.
\end{lem}

\begin{proof}
Let $M_t$ be the polydule given by $M_t(S,i)=t^{|S|-1}M(S,i)$ for
all $|S|>1$. We have a strict isomorphism $M_t\cong M$ for all
$t\neq 0$ and $M_0=H^*(M)$. Therefore $M\dege H^*(M)$.
\end{proof}

Assuming for a moment the results of the next section, we now show
that degeneration of polyrepresentations is equivalent to
degeneration of polydules.

\begin{lem} \label{replemma}
Let $M,N\in \Rep_\infty(\Lambda,\underline{e})$. Then $M\dege N$ if
and only if $N\in \overline{\mathcal{G}_Q\cdot M}$.
\end{lem}

\begin{proof}
Assume $N\in \overline{\mathcal{G}_Q\cdot M}$. We have
$\mathcal{G}_Q\cdot M \subseteq \mathcal{G}\cdot M$ and therefore
$\overline{\mathcal{G}_Q\cdot M} \subseteq \overline{\mathcal{G}\cdot
  M}$. Hence $M\dege N$.

Now assume that $M\dege N$. So there exists a polydule $Z$ in
$\mathcal D^f_\infty(\L)$ and an exact triangle $$Z \longrightarrow
M\oplus Z \longrightarrow N \longrightarrow Z[1]$$ by Theorem
\ref{zwar}. We may assume that $Z$ is a polyrepresentation. Now, by
the proof of Theorem \ref{riedt} there exists a family of
representations $M_t$ in $\overline{\mathcal{G} _Q\cdot M}$, such
that $M_t\cong M$ for all except a finite number of $t$, and
$M_0\cong N$. This shows that $N\in \overline{\mathcal{G}_Q\cdot
M}$.
\end{proof}

\subsection{Polydules with fixed homology}\label{fix}

Let $X \in \mod (\Lambda,\underline{d})$ be a graded module with
graded dimension $\underline{d}$. Let $\mathcal{M}_X=\pi^{-1}(X)$ be
the variety of polydules with fixed homology $X$. Let
$\mathcal{G}_X$ be the setwise stabilizer of $\mathcal{M}_X$,
consisting of group elements with $F(\empt,i)\in
\Aut_\Lambda(H^i(X))$. There is a bijection between the orbits under
the induced action of $\mathcal{G}_X$ on $\mathcal{M}_X$ and the
isomorphism classes of polydules with homology equal to $X$. The
variety $\mathcal{M}_X$ is connected with unique closed orbit the
orbit of $X$. The following lemma may be proven similar to Lemma
\ref{replemma} above. We skip the details.

\begin{lem}\label{homlemma}
Let $M,N\in \mathcal{M}_X$. Then $M\dege N$ if and only if $N\in
\overline{\mathcal{G}_X\cdot M}$.
\end{lem}

As a special case we can consider polydules with fixed homology in two
degrees. That is, $X=M[a]\oplus N[b]$ for two $\Lambda$-modules $M$
and $N$ and degrees $a$ and $b$. Then the equations defining
$\mathcal{M}_X$ are affine equations, and therefore $\mathcal{M}_X$ is
isomorphic to an affine space. The following lemma is a useful application
of varieties of polydules to the problem of checking the vanishing of
extension spaces between $\Lambda$-modules.

\begin{lem} \label{extlemma}
Let $M$ and $N$ be two finitely generated $\Lambda$-modules and let
$X=M\oplus N[i]$, where $i \geq 0$. Then $\Ext^{i+1}_\L(M,N)=0$ if
and only if $\mathcal{G}_X\cdot X = \mathcal{M}_X$.
\end{lem}

\begin{proof}
We have $\Ext^{i+1}_\L (M,N)=\Hom(M,N[i+1]) \simeq
\Hom(M[-1],N[i])$. If $f \colon M[-1] \pil N[i]$ is a morphism, then
the cone of $f$ is in $\mathcal{M}_X$, and every point in
$\mathcal{M}_X$ can be obtained as the cone of such a morphism. The
cone of $f$ is quasi-isomorphic to $X$ if and only if $f$ is
null-homotopic. The lemma follows.
\end{proof}

Suppose $\Lambda=kQ/I$ for a quiver $Q=(Q_0,Q_1)$ and an admissible
ideal $I\subseteq kQ$. We can combine the constructions in this
section and consider polyrepresentations with fixed homology. Let $X
\in \Rep(\L,\underline{e})$ be a graded representation with sequence
of dimension vectors
$\underline{e}=(\underline{e}^a,\ldots,\underline{e}^b)$. Then
define $$\mathcal M^Q_X=\mathcal M_X \cap
\Rep_{\infty}(\L,\underline{e}),$$ or equivalently $\mathcal
M^Q_X=(\pi')^{-1}(X)$, where $\pi'$ is the restricted function $\pi
'=\pi |_{\Rep_{\infty}(\L,\underline{e})} \colon
\Rep_{\infty}(\L,\underline{e}) \pil \Rep(\L,\underline{e})$. There
is a group $\mathcal G^Q_X =\mathcal G_X \cap \mathcal
G_{Q,\underline e}$ acting on $\mathcal M^Q_X$ and orbits correspond
to isomorphism classes of polyrepresentations. Analogues of Lemmas
\ref{homlemma} and \ref{extlemma} hold in this situation.

\section{Characterisation of orbit closure}\label{orbit}

In this section we prove our main result, an algebraic
characterisation of degeneration (orbit closure) in the variety
$\mathcal M_{\L,\underline d}$. Keeping in mind the equivalence $\Db
(\mod \L) \xpil \sim \mathcal D^f_\infty(\L)$, this should be seen
as a derived category version of a celebrated theorem by Zwara
\cite{Zwa} and Riedtmann \cite{Riedt}. In fact the Riedtmann-Zwara
Theorem is a special case of the result we prove here.

In \cite{Orb} the authors gave a very general theorem characterising
orbit closures in varieties. It takes the following form when
applied to our situation.

\begin{thm}\label{dvr}
Let $M, N \in \mathcal M_{\L,\underline d}$. Then $M \dege N$ if and
only if there is a discrete valuation $k$-algebra $R$ (residue field
$k$, quotient field $K$ finitely generated and of transcendence
degree one over $k$) and a $\L \otimes_k R$-polydule structure with
$m_1=0$ on the graded $R$-module $Q= \oplus_{i=a}^b R^{d_i}$ such
that $Q \otimes_R k = N$ (via the identification $R \otimes_R k=k$)
and $Q \otimes_R K \simeq M \otimes_k K$ as $\L \otimes_k
K$-polydules.
\end{thm}

\begin{proof}
The statement follows from \cite[1.2]{Orb}.
\end{proof}

Using this as a starting point, we can prove the "Zwara direction"
of our characterisation. Zwara \cite{Zwa} was the first to prove the
analogous statement for module varieties. Yoshino \cite{Yos} defined
module degeneration for any $k$-algebra, also when there is no
underlying variety, and showed that Zwara's result was still valid
in this generality \cite[2.2]{Yos}. Our proof, although completely
different due to the triangulated setting, borrows some ideas from
\cite{Yos}.

\begin{thm}["Zwara direction"] \label{zwar}
Let $M, N \in \mathcal M_{\L,\underline d}$. If $M \dege N$, then
there is an object $Z$ in $\mathcal D ^f_\infty(\L)$ and an exact
triangle
$$Z \pil M \oplus Z \pil N \pil Z[1].$$
\end{thm}

\begin{proof}
Suppose $M \dege N$. Let $R$ be a discrete valuation $k$-algebra as
in Theorem \ref{dvr} and let $Q= \oplus_{i=a}^b R^{d_i}$ be a $\L
\otimes_k R$-polydule such that $m_1^Q=0$, $Q \otimes_R k = N$ and
$Q \otimes_R K \simeq M \otimes_k K$ as $\L \otimes_k K$-polydules.
Let $t$ denote a generator of the maximal ideal in $R$.

Let $i \colon M \otimes_k R \pil M \otimes_k K$ be the strict
inclusion of $\L \otimes_k R$-polydules and let $f$ be the
composition
$$M \otimes_k R \xpil i M \otimes_k K \xpil \sim Q \otimes_R K.$$
Since $M \otimes_k R$ and $Q \otimes_R K$ have support in only a
finite number of degrees, we know that $f_i=0$ for $i>>0$. Since in
addition $M \otimes_k R$ is finitely generated as an $R$-module,
there is an $n$ such that the morphism $f$ factors through the
strict inclusion $Q \otimes_R (1/t^n) R \pil Q \otimes_R K$. There
is a strict isomorphism $Q \otimes_R (1/t^n) R \simeq Q$. We have
the following morphism of triangles in $\mathcal D_\infty(\L)$.
$$\xymatrix{M \otimes_k R \ar[r]^g \ar@{=}[d] & Q \ar[r] \ar[d]
& \cone g \ar[r]^(.4)v \ar[d] & (M \otimes_k R)[1] \ar@{=}[d] \\
M \otimes_k R \ar[r]^f & Q \otimes_R K \ar[r] & \cone f \ar[r]^(.4)u
& (M \otimes_k R)[1]}$$ Since $f$ is a section, we have $u=0$ by
\cite[I.4.1]{Hap}. Since $u=0$, also $v=0$ and $g$ is a section. It
follows that $g_1 \colon M \otimes_k R \pil Q$ is injective. Since
$M \otimes_k R$ and $Q$ have the same rank over $R$, it follows by
the construction of the cone that $\cone g$ has finite dimensional
total homology.

Now consider the strict inclusion of $\L$-polydules $j \colon M
\otimes_k tR \pil  M \otimes_k R$. From the exact sequence of vector
spaces $0 \pil tR \pil R \pil k \pil 0$, we get an exact sequence of
$\L$-polydules $0 \pil M \otimes_k tR \xpil j M \otimes_k R \pil M
\pil 0$, so $M \simeq \cone j$ in $\mathcal D_\infty(\L)$. By the
octahedral axiom, we get the following diagram of exact (marked by
$\triangle$) and commuting triangles in $\mathcal D_\infty(\L)$.

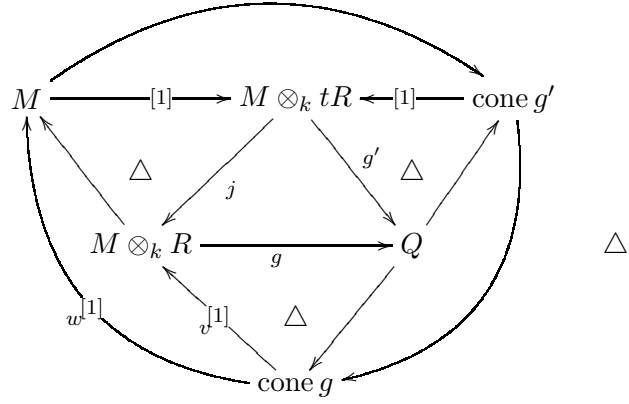
$$\xymatrix@=0.4cm{&&  \\ M \ar[rr]|{[1]} \ar@/^3pc/[rrrr] && M \otimes_k tR
\ar[ddl]^j
 \ar[ddr]^{g'}    &&  \cone  g'
\ar[ll]|{[1]} \ar@/^3pc/[ddddll]\\ &\triangle  && \triangle \\& M
\otimes_k R \ar[rr]_g \ar[uul] &  & Q \ar[uur] \ar[ddl] && \triangle
\\ &  & \triangle
\\ && \cone g \ar[uul] | {[1]}^v \ar@/^3pc/[uuuull]| {[1]}^w}$$

So there is an exact triangle $M \pil \cone g' \pil \cone g \xpil w
M[1]$. Since $v=0$, also $w=0$ and $\cone g' \simeq M \oplus \cone
g$.

From the commutative square $$\xymatrix{M \otimes_k R \ar[r]^(.6)g
\ar[d]^{\tilde t} & Q \ar[d]^{\tilde t}\\
M \otimes_k tR \ar[r]^(.6){g'} & Q},$$ where $\tilde t$ denotes
multiplication by $t$, we get the commutative diagram
$$\xymatrix{M \otimes_k R \ar[r]^g \ar[d]^{\tilde t} & Q \ar[r]
\ar[d]^{\tilde t}
& \cone g \ar[r]^(.4)v \ar[d]^h & (M \otimes_k R)[1] \\
M \otimes_k tR \ar[r]^{g'} \ar[d] & Q \ar[r] \ar[d] & \cone g'
\ar[r] \ar[d] &
(M \otimes_k tR)[1]\\0 \ar[r] \ar[d] & Q \otimes_R k \ar[r] \ar[d] & \cone h \ar[r] \ar[d] & 0\\
(M \otimes_k R) [1] & Q [1] & (\cone g)[1] }$$ where the rows and
columns are exact triangles \cite[Ex.10.2.6.]{Wei}. We see that
$\cone h \simeq Q \otimes_R k = N$, and there is an exact triangle
$$\cone g \xpil h \cone g' \pil \cone h \pil (\cone g)[1]$$ in
$\mathcal D_\infty(\L)$. We put $Z=\cone g$, which is an object in
$\mathcal D^f_{\infty}(\L)$. The triangle above is isomorphic to $Z
\pil M \oplus Z \pil N \pil Z[1]$.
\end{proof}

In order to prove the converse of this theorem, we first make the
useful observation that the existence of an exact triangle of the
form in question implies the existence of a certain exact sequence
of polydules.

\begin{lem}
Let $M,N \in \mathcal M_{\L,\underline d}$. Suppose there is an
exact triangle of the form $Z \pil M \oplus Z \pil N \pil Z[1]$ in
$\mathcal D^f_\infty(\L)$.

Then there is an exact sequence of finite dimensional $\L$-polydules
of the form $0 \pil Z \xpil f M' \oplus Z' \xpil g N \pil 0$ with
$m_1^{Z}=m_1^{Z'}=m_1^{M'}=0$ and where $Z \simeq Z'$ and $M \simeq
M'$. (These isomorphisms are possibly non-strict.)
\end{lem}

\begin{proof}
Suppose there is a triangle $Z \pil  M \oplus Z \pil N \xpil h Z[1]$
in $\mathcal D^f_\infty(\L)$. We may assume that $Z$ is represented
in such a way that $m_1^Z=0$. The morphism $h[-1] \colon N[-1] \pil
Z$ gives rise to an exact sequence $0 \pil Z \xpil f \cone (h[-1])
\xpil g N \pil 0$. Since $\dim_k M^i=\dim_k N^i$ for all $i \in \Z$,
also $\dim_k (M \oplus Z)^i=\dim_k (\cone (h[-1]))^i$ for all $i \in
\Z$. The quasi-isomorphism $M \oplus Z \pil \cone (h[-1])$ must
therefore be a (possibly non-strict) isomorphism and $m_1^{\cone
(h[-1])}=0$.
\end{proof}

We are now ready to prove the converse of Theorem \ref{zwar}. Our
proof is similar to Riedtmann's original proof for module varieties
\cite[3.4]{Riedt}.

\begin{thm}["Riedtmann direction"]\label{riedt}
Let $M,N \in \mathcal M_{\L,\underline d}$. If there is an exact
triangle $Z \pil M \oplus Z \pil N \pil Z[1]$ in $\mathcal
D^f_\infty(\L)$, then $M \dege N$.
\end{thm}

\begin{proof}
Suppose there is a triangle $Z \pil  M \oplus Z \pil N \xpil h Z[1]$
in $\mathcal D^f_\infty(\L)$. From the previous lemma it follows
that there is an exact sequence of finite dimensional $\L$-polydules
of the form $0 \pil Z \xpil f M' \oplus Z' \xpil g N \pil 0$ with
$m_1^{Z}=m_1^{Z'}=m_1^{M'}=0$ and where $v \colon Z \simeq Z'$ and
$M \simeq M'$. We now show that this implies $M \dege N$.

Let $W= \Im f$. For each $i \in Z$, choose a $k$-complement $C^i$ of
$W^i$ in $M^i \oplus Z^i$, so $W^i \oplus C^i \simeq (M')^i \oplus
(Z')^i$ as $k$-vector spaces. Let $\mathcal F$ be the set of all
strict morphisms $\tilde f \colon Z \pil M' \oplus Z'$ such that
$\tilde f$ is injective and $(\Im \tilde f) \cap C = 0$. These are
open conditions, so $\mathcal F$ is a variety.

For each $\tilde f \in \mathcal F$ we define an $\ainf$-module
structure on $C$ and a strict morphism $\tilde g \colon M' \oplus Z'
\pil C$ such that $0 \pil Z \xpil {\tilde f} M' \oplus Z' \xpil
{\tilde g} C \pil 0$ becomes an exact sequence. Let $\tilde f=
\left ( \begin{smallmatrix} \tilde \phi \\
\tilde \psi
\end{smallmatrix} \right ) \colon Z \pil W \oplus C$. We assume $\tilde f \in \mathcal
F$, so $\tilde \phi$ is an invertible map of graded vector spaces.
We define $$\tilde g=\left ( \begin{matrix} \tilde \psi (\tilde
\phi)^{-1} & \mathtt 1_C
\end{matrix} \right ) \colon W \oplus C \pil C$$ and $$m_n^C(\lambda_{n-1},
\ldots,\lambda_1,a)= \tilde g (m_n^{M' \oplus Z'}(\lambda_{n-1},
\ldots,\lambda_1,a))$$ for all $n \geq 2$, $a \in C$ and $\lambda_1,
\ldots,\lambda_{n-1} \in \L$. Since $\dim_k C^i=\dim_k N^i$ for all
$i \in \Z$, this determines a morphism of varieties $\theta \colon
\mathcal F \pil \mathcal M_{\L,\underline d}$. This morphism is
defined in such a way that $\theta(\tilde f)$ corresponds to a
polydule quasi-isomorphic to $\cone \tilde f$.

There is an open set $U \in \mathbb A_1$ containing $0$ such that
the morphism $f_{(t)}=f+t \cdot {\left ( \begin{smallmatrix}  0 \\
v
\end{smallmatrix} \right ) } \colon Z \pil M' \oplus Z'$ belongs to $\mathcal F$ whenever $t
\in U$. Except for a finite number of points in $U$, the second
component of $f_{(t)}$ has full rank and $\cone f_{(t)} \simeq M'
\simeq M$. If $t=0$, then $\cone f_{(0)}= \cone f \simeq N$. So $M
\dege N$.

\end{proof}

Combining the results so far in this section we get exactly the
algebraic characterisation of degeneration we were aiming for.

\begin{thm}\label{rizw}
Let $M,N \in \mathcal M_{\L,\underline d}$. Then the following are
equivalent.
\begin{itemize}
\item[(a)]$M \dege N$.
\item[(b)]There is an object $Z$ in $\mathcal
D^f_\infty(\L)$ and an exact triangle
$$Z \pil M \oplus Z \pil N \pil Z[1].$$
\item[(c)] There is an exact sequence of finite dimensional
$\L$-polydules of the form $$0 \pil Z
\pil  M' \oplus Z' \pil  N \pil 0$$ with
$m_1^{Z}=m_1^{Z'}=m_1^{M'}=0$ and where $Z \simeq Z'$ and $M \simeq
M'$.
\end{itemize}
\end{thm}

If $a=b=0$ and $\underline d=(d)$, then $\mathcal M_{\L,\underline
d}$ is equal to the module variety $\mod(\L,d)$. We get the
Riedtmann-Zwara Theorem (\cite{Zwa}, \cite{Riedt}) as a corollary.

\begin{cor}[Riedtmann-Zwara Theorem]
Let $\L$ be a finite dimensional $k$-algebra and let $M$ and $N$ be
two $d$-dimensional $\L$-modules. Then $M \dege N$ if and only if
there exists a finite dimensional $\L$-module $Y$ and a short-exact
sequence $0 \pil Y \pil M \oplus Y \pil N \pil 0$.
\end{cor}

\begin{proof}
If $v \colon Z \pil Z'$ is an isomorphism of polydules, then it
induces a bijective $\L$-morphism (so a $\L$-isomorphism) $v_1
\colon Z^0 \pil (Z')^0$. From the equivalence between (a) and (c) in
the theorem it follows that $M \dege N$ if and only if there exists
a finite dimensional $\L$-module $Z^0$ and an exact sequence of
$\L$-modules $0 \pil Z^0 \pil M \oplus Z^0 \pil N \pil 0$.
\end{proof}

In \cite{Der}, another theory of degeneration in derived categories
was developed. There degenerations are studied in the space
$\comj^{\underline d}$, a space which parameterise right bounded
complexes of projective modules with so-called dimension array
$\underline d$. There is a group $\mathsf G$ acting on this space
and $\mathsf G$-orbits correspond to isomorphism classes in
$\mathcal D (\mod \L)$ where the dimension array appears (see
\cite{Der} for details). One then defines $M \leq_{\top} N$ if $N$
is in the closure of the $\mathsf G$-orbit of $M$. This concept of
degeneration is compatible with the one in the present paper.

\begin{cor}
Let $X$ and $Y$ be two right bounded complexes of finitely generated
projective $\L$-modules with same dimension array $\underline d$.
Suppose $\dim_k H^i X=\dim_k H^i Y$ for all $i \in \Z$ and that
$\sum_{i \in \Z} \dim_k H^i X < \infty$. Let $M_X$ and $M_Y$ be
$\L$-polydules with $m_1^{M_X}=0= m_1^{M_Y}$ which are
quasi-isomorphic to $X$ and $Y$ respectively.

Then $M_X \dege M_Y$ if and only if $X \leq_{\top} Y$.
\end{cor}

\begin{proof}
The relation $X \leq_{\top} Y$ can be characterised by the existence
of a right bounded complex $Z$ in $\mathcal D (\mod \L)$ and an
exact triangle $Z \pil X \oplus Z \pil Y \pil Z[1]$ (see
\cite{Der}). Since both $X$ and $Y$ have bounded homology, also $Z$
can be chosen in $\Db(\mod \L)$. Now the statement follows from
Theorem \ref{rizw} and the equivalence of triangulated categories
$\Db (\mod \L) \xpil \sim \mathcal D^f_\infty(\L)$.
\end{proof}

\section{Examples}\label{ex}

We now give some examples illustrating the main theorem and also the
different varieties constructed in section \ref{simp}.

\subsection{Example 1}

Let $M \in \mathcal M_{\L,\underline d}$. From Lemma \ref{deglemma}
we know that $M\dege H^*M$. We find a polydule $Z$ and a triangle
$Z\longrightarrow M\oplus Z \longrightarrow H^*M \longrightarrow
Z[1]$ in this case. Assume that $H^i M$ is non-zero only in an
interval $s \leq i \leq t$. We use the notation $Z_i= \tau^{\leq
i}M$ for the truncated polydules. (Truncation is straightforward
when $m_1^M=0$.) We have an exact triangle
$$\eta_{t-1} \colon Z_{t-1}\longrightarrow M \longrightarrow H^tM
\longrightarrow Z_{t-1}[1]$$ and by induction exact triangles
$$\eta_{i} \colon Z_{i}\longrightarrow Z_{i+1} \longrightarrow H^{i+1}M
\longrightarrow Z_{i}[1]$$ for $s-1 \leq i <t$. Then by taking the
direct sum $\eta=\eta_{t-1}\oplus \dots \oplus \eta_{s-1}$ of exact
triangles and letting $Z=Z_{t-1}\oplus Z_{t-2} \oplus \dots \oplus
Z_s \oplus Z_{s-1}$, where $Z_{s-1}=0$, we have the exact triangle
$$\eta \colon Z\longrightarrow M\oplus Z \longrightarrow H^*M
\longrightarrow Z[1].$$

\subsection{Example 2}

Let $A$ be the algebra given by the quiver
$$\xymatrix{1 \ar[r]^{\alpha} & 2 \ar[r]^{\beta} & 3 \ar[r]^{\gamma} &
4 \ar[r]^{\delta} & 5 }$$ with relations $\beta \alpha=0= \delta
\gamma$. The global dimension of $A$ is $2$. Let
$$X_1=S_5[2]\oplus S_3[1]\oplus S_1,$$ where $S_i$ denotes the
simple module corresponding to vertex $i$. We compute the variety
$\mathcal{M}_{X_1}^Q$ of polyrepresentations with homology module
$X_1$, as defined in subsection \ref{fix}. Let
\begin{align*}
x &=M((\beta,\alpha),0),\\
y &=M((\delta,\gamma),-1),\\
z &=M((\delta,\gamma\beta,\alpha),0).
\end{align*}
(We identify $M(S,i)$ with its unique non-zero block, a natural
thing to do when working with polyrepresentations.) Then
$\mathcal{M}_{X_1}^Q$ is the subvariety of $\mathbb{A}^3$ given by
the equation $z=yx$. The group
$$\mathcal{G}_{X_1}^Q=(k^*)^3$$ acts by
$(a,b,c)(x,y,z)=(bxc^{-1},ayb^{-1},azc^{-1})$. There is only the
action given by change of basis in this example. There are four
orbits given by the polydules
$$M_1 \colon x=1, y=1, z=1,$$
$$M_2 \colon x=1, y=0, z=0,$$
$$M_3 \colon x=0, y=1, z=0,$$
$$M_4 \colon x=0, y=0, z=0.$$ The
Hasse diagram of degeneration is given by $$
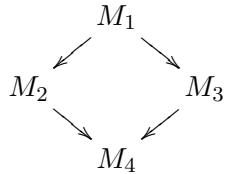
\xymatrix@=0.4cm{& M_1
\ar[dr] \ar[dl] & \\ M_2\ar[dr] && M_3
  \ar[dl] \\ & M_4 &}$$

Consider the degeneration $M_1 \dege M_3$. Let $\L= \tau^{\leq -1}
M_3$. We have a decomposition $M_3 = S_1 \oplus L$. There is an
exact sequence of polydules $0 \pil L \xpil f M_1 \xpil g S_1 \pil
0$. Using a standard trick, we get an exact sequence of polydules
$$0 \pil
 L \xpil {\left( \begin{smallmatrix} f & 0
\end{smallmatrix} \right )} M_1 \oplus L \xpil h M_3 \pil 0$$
where $h=\left ( \begin{matrix} g & 0 \\ 0 & \mathtt 1_L
\end{matrix}\right ) \colon M_1\oplus L \longrightarrow  S_1\oplus
L$. This is a sequence of the same form as in Theorem \ref{rizw}.

\subsection{Example 3}

Let $A$ be as in the previous example. Let
$\underline{e}^{-1}=(1,0,0,0,0)$ and $\underline{e}^0=(0,0,1,1,0)$
and be two dimension vectors and let
$\underline{e}=(\underline{e}^{-1},\underline{e}^0)$. We consider
the variety $\Rep_\infty(A,\underline{e})$ defined in subsection
\ref{quiv}. Let
\begin{align*}
x &=M((\beta,\alpha),0),\\
y &=M((\gamma\beta,\alpha),0),\\
z &=M((\gamma),-1).
\end{align*}
Then $\Rep_\infty(A,\underline{e})$ is given by the triples
$(x,y,z)$ with $y=-zx$. The group $G_{Q,\underline{e}}=(k^*)^3$ acts
by $(a,b,c)(x,y,z)=(bxa^{-1},cya^{-1},czb^{-1})$. There is only the
action given by change of basis in this example. There are four
orbits given by the polydules
\begin{align*}
M_1 &\colon x=1, y=-1, z=1,\\
M_2 &\colon x=1, y=0, z=0,\\
M_3 &\colon x=0, y=0, z=1,\\
M_4 &\colon x=0, y=0, z=0.
\end{align*}
The Hasse diagram of degeneration is given by $$
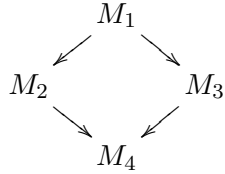
\xymatrix@=0.4cm{&
M_1 \ar[dr] \ar[dl] & \\ M_2\ar[dr] && M_3
  \ar[dl] \\ & M_4 &}$$

Let $M$ be the representation $$\xymatrix{0
  \ar[r] & 0 \ar[r] & k \ar[r]^1 & k \ar[r] & 0}.$$ Then
$$H^*(M_1)=S_1\oplus M[1] = H^*(M_3)$$ and $$H^*(M_2)=H^*(M_4)=S_1\oplus
(S_3\oplus S_4)[1].$$

\subsection{Example 4}

Let $A$ be the algebra of the previous two examples. Let $X_2 =
S_j[2] \oplus S_1$ for some $j=1,2,3,4,5$. Since the global
dimension of $A$ is $2$, we see that $\Ext_A^3(S_1,S_j)=0$. So
according to Lemma \ref{extlemma} (or more precisely its version for
polyrepresentations) we have $\mathcal{G}_{X_2}^Q \cdot X_2 =
\mathcal M_{X_2}^Q$. We show this in detail. Since the
multiplication $m_4$ takes three arguments from $A$, we must have
$j=4$ or $5$ in order to have a non-trivial variety
$\mathcal{M}_{X_2}^Q$. If $j=5$, then $\mathcal{M}_{X_2}^Q$ is given
by $x=M((\delta,\gamma\beta,\alpha),0)$ with the equation $x=0$,
which gives the trivial variety. If $j=4$, then
$\mathcal{M}_{X_2}^Q$ is given by the affine line
$x=M((\gamma,\beta,\alpha),0)$ with no relations. The group
$\mathcal{G}_{X_2}^Q=(k^*)^2\times k$ acts by
$$(a,b,c)\cdot t = (ax+c)b^{-1}.$$ Hence there is only one orbit.
This shows that $\mathcal{G}_{X_2}^Q\cdot X_2 = \mathcal M_{X_2}^Q$.


\begin{thebibliography}{ABCD}

\bibitem[ARS] {Ars} M. Auslander, I. Reiten, S. O. Smal\o,
  \emph{Representation theory of Artin algebras.}
Cambridge Studies in Advanced Mathematics, 36. Cambridge University
Press, Cambridge, 1997. xiv+425 pp.

\bibitem[GO'H] {Orb} F. Gr\" unewald, J. O'Halloran, \emph{A
characterization of orbit closure and applications}, J. Algebra 116
(1988), no. 1, 163--175.

\bibitem[H] {Hap} D. Happel, \emph{Triangulated categories in the
 representation theory of finite-dimensional algebras.}
London Mathematical Society Lecture Note Series, 119. Cambridge
University Press, Cambridge, 1988. x+208 pp.

\bibitem[JSZ] {Der} B. T. Jensen, X. Su, A. Zimmermann,
\emph{Degenerations for derived categories}, J. Pure Appl. Algebra
198 (2005), no. 1-3, 281--295.

\bibitem[K] {Ainf} B. Keller, \emph{Introduction to $A$-infinity
algebras and modules}, Homology Homotopy Appl. 3 (2001), no. 1,
1--35.

\bibitem[L-H] {Doc} K. Lef\`{e}vre-Hasegawa, \emph{Sur les
$\ainf$-cat\'{e}gories}, Th\`{e}se de Doctorat, Universite Paris 7,
2003.

\bibitem[R] {Riedt} C. Riedtmann,
\emph{Degenerations for representations of quivers with relations.}
Ann. Sci. École Norm. Sup. (4) 19 (1986), no. 2, 275--301.

\bibitem[W] {Wei} C. A. Weibel, \emph{An introduction to homological
algebra.} Cambridge Studies in Advanced Mathematics, 38. Cambridge
University Press, Cambridge, 1994. xiv+450 pp.

\bibitem[Y] {Yos} Y. Yoshino, \emph{On degeneration of modules},
J. Algebra 278 (2004), no. 1, 217--226.

\bibitem[Z] {Zwa} G. Zwara, \emph{Degenerations of finite
dimensional modules are given by extensions}, Compositio Math. 121
(2000), no. 2, 205--218.



\end{thebibliography}
\end{document}